\documentclass[a4paper,12pt]{amsart}
\usepackage{graphicx}
\usepackage{latexsym,amssymb,amsmath,mathrsfs}
\usepackage{amsthm}
\usepackage{mathrsfs}
\usepackage{color}
\usepackage[english]{babel}
\usepackage{comment}
\usepackage[top=30truemm,bottom=30truemm,left=20truemm,right=20truemm]{geometry}
\usepackage[hidelinks]{hyperref}
\usepackage[abbrev]{amsrefs}
\newtheorem{theorem}{Theorem}[section]
\newtheorem{lemma}[theorem]{Lemma}
\newtheorem{proposition}[theorem]{Proposition}

\theoremstyle{remark}
\newtheorem{remark}[theorem]{Remark}

\theoremstyle{remark}

\DeclareMathOperator{\Ric}{Ric}

\newcommand{\Hess}{\mathrm{Hess}\,}
\renewcommand{\d}{\mathrm{d}}

\title[A remark on the compactness]{A remark on the compactness of manifolds \\with non-negative Ricci curvature\\ and uniformly convex boundary}
\author{Yasuaki Fujitani}
\address{Graduate School of Mathematical Sciences, The University of Tokyo, Komaba, Tokyo, 153-8914, Japan} 
\email{yasuakifujitani@g.ecc.u-tokyo.ac.jp}
\begin{document}
\begin{abstract}
Yan--Zhu \cite{yan2026nonnegative} showed the compactness of complete Riemannian manifolds with non-negative Ricci curvature and uniformly convex boundary.
In this paper,
we give an alternative proof of this compactness theorem.
\end{abstract}
\maketitle
\section{Introduction}
Throughout this paper,
$(M,g)$ denotes an $n$-dimensional complete connected Riemannian manifold with smooth nonempty boundary and $n \geq 2$.
We denote the induced metric and the second fundamental form on $\partial M$ by $g_{\partial M}$ and $\mathrm{II}_{g,\partial M}$,
respectively.
The following compactness theorem,
conjectured by Li \cite[Conjecture 1.3]{li2014sharp},
was proved by Yan--Zhu \cite[Theorem 1.2]{yan2026nonnegative}:
\begin{theorem}[\cite{yan2026nonnegative}]\label{thm:yan-zhu}
Let $(M,g)$ be an $n$-dimensional complete Riemannian manifold with boundary.
For $\kappa > 0$,
we assume
\begin{align}\label{eq:curvature}
    \Ric_g \geq 0,\quad \mathrm{II}_{g,\partial M} \geq \kappa g_{\partial M}.
\end{align}
Then $M$ is compact.
\end{theorem}
Although the proof of Yan--Zhu \cite{yan2026nonnegative} uses monotone quantities and the $L^2$-trace Poincar\'{e} inequality in \cite[Theorem A.1]{yan2026nonnegative},
we give a different proof by comparing lower and upper bounds for the infimum of the Neumann Rayleigh quotient over nonzero Sobolev functions orthogonal to the kernel of the Neumann Laplacian.
We denote this infimum of the Neumann Rayleigh quotient by $\lambda_1^{\mathrm{Neu}}(M,g)$.
We refer to Section \ref{sec:proof} for its precise definition.
We obtain the lower bound of $\lambda_1^{\mathrm{Neu}}(M,g)$ by using the $L^1$-trace Poincar\'{e} inequality in Yan--Zhu \cite[Proposition A.3]{yan2026nonnegative}.
Its constant-function case recovers the Ros inequality \cite{ros1987compact}.
We note that the proof of Theorem \ref{thm:yan-zhu} in \cite{yan2026nonnegative} does not use \cite[Proposition A.3]{yan2026nonnegative}.

The quantity $\lambda_1^{\mathrm{Neu}}(M,g)$ coincides with the first nonzero Neumann eigenvalue if $M$ is compact.
The Reilly formula is a classical tool to obtain lower bounds for Dirichlet and Neumann Laplace eigenvalues on compact manifolds with boundary.
For example,
we refer to Reilly \cite{reilly1977applications} for Dirichlet eigenvalue estimates,
and we refer to Escobar \cite{escobar1990uniqueness} and Xia \cite{xia1991first} for Neumann eigenvalue estimates.
Different from our setting,
their results assume positive Ricci curvature.

Rayleigh-type variational quantities have also long been studied on noncompact manifolds.
For manifolds without boundary,
McKean \cite{mckean1970upper} and Cheng \cite{cheng1975eigenvalue} obtained lower and upper bounds,
respectively,
under some curvature assumptions.
For noncompact manifolds with boundary,
in a more general weighted setting,
Wang \cite{wang1999existence} investigated lower and upper bounds for such a variational quantity (see also Chen--Wang \cite{chen1997general}).
In their setting,
the Riemannian volume measure is replaced by a more general weighted measure,
and the Laplace operator is replaced by a symmetric diffusion operator (see also Bakry--\'{E}mery \cite{bakry1985diffusions}).
In the unweighted finite volume case,
the variational quantity in \cite{wang1999existence} coincides with $\lambda_1^{\mathrm{Neu}}(M,g)$ (see also Remark \ref{rem:lambda} below).
Their arguments use methods from probability theory.
Our lower and upper bounds for $\lambda_1^{\mathrm{Neu}}(M,g)$ are different from those results mentioned above.


In the definition of $\lambda_1^{\mathrm{Neu}}(M,g)$,
the Neumann Laplacian is used.
The functional-analytic properties of the Laplacian on non-compact manifolds such as essential self-adjointness have also been investigated.
For manifolds without boundary,
the theory of essential self-adjointness goes back to the work of Gaffney \cite{gaffney1951harmonic,gaffney1954special,gaffney1955hilbert} and Roelcke \cite{roelcke1960laplace} (see also Shubin \cite{shubin1992spectral}).
Furthermore,
Strichartz \cite{strichartz1983analysis} gave a new proof of this essential self-adjointness.
For manifolds with boundary,
essential self-adjointness was proved by Bianchi--G\"{u}neysu--Setti \cite{bianchi2025neumann},
which we use to obtain the lower bound for $\lambda_1^{\mathrm{Neu}}(M,g)$.

To obtain the lower bound for $\lambda_1^{\mathrm{Neu}}(M,g)$,
we also use the Reilly formula \cite{reilly1977applications},
the $L^1$-trace Poincar\'{e} inequality in Yan--Zhu \cite[Proposition A.3]{yan2026nonnegative},
and the inscribed radius estimate (see Kasue \cite{kasue1983ricci} and Li \cite{li2014sharp}).
Indeed,
we first obtain an estimate for smooth functions by using these geometric inequalities.
After that,
together with the approximation property by the self-adjointness theorem in Bianchi--G\"{u}neysu--Setti \cite{bianchi2025neumann},
we obtain our lower bound estimate for a more general class of functions.

To obtain the upper bound for $\lambda_1^{\mathrm{Neu}}(M,g)$,
we use cut-off functions constructed from the distance function and Bishop--Gromov volume comparison.
We note that there is also a proof of a Cheng-type upper bound \cite{cheng1975eigenvalue} via the volume comparison (see e.g., Brooks \cite{brooks1981relation} and Gallot--Hulin--Lafontaine \cite[Section 4.G]{gallot2004riemannian}). 
Our argument for the upper bound is similar to that proof using the volume comparison.

The organization of this paper is as follows.
In Section \ref{sec:Preliminary},
we introduce terminology and properties used in Section \ref{sec:proof}.
Subsection \ref{subsec:functional-analysis} introduces terminology in functional analysis,
and Subsection \ref{subsec:functional-analysis-boundary} introduces the essential self-adjointness theorem.
Moreover,
Subsection \ref{subsec:geometric-inequality} introduces Bishop--Gromov volume comparison,
the inscribed radius estimate,
and the $L^1$-trace Poincar\'{e} inequality.
In Section \ref{sec:proof},
we obtain lower and upper bounds for $\lambda_1^{\mathrm{Neu}}(M,g)$,
and combine them to give the alternative proof of Theorem \ref{thm:yan-zhu}.
\section{Preliminaries}\label{sec:Preliminary}
In this section,
we collect the analytic and geometric ingredients for Section \ref{sec:proof}.
\subsection{Functional analysis}\label{subsec:functional-analysis}
In this subsection,
we introduce terminology of functional analysis \textcolor{black}{(see e.g., Brezis \cite{brezis2011functional} and Schm\"{u}dgen \cite{schmudgen2012unbounded})}.
Let $\mathcal{H}$ be a real Hilbert space,
and let $\langle\cdot, \cdot\rangle_\mathcal{H}$ be the inner product of $\mathcal{H}$.
Let 
\begin{align*}
    A : \mathrm{Dom}(A)  \longrightarrow \mathcal{H}
\end{align*}
be a \textcolor{black}{densely} defined linear operator.
Here,
$\mathrm{Dom}(A)$ is a dense linear subspace of $\mathcal{H}$.
For a linear operator $B$,
we write $A \subset B$ if $\mathrm{Dom}(A) \subset \mathrm{Dom}(B)$ and $A x = B x$ for every $x \in \mathrm{Dom}(A)$.
In this case,
$B$ is called an \textit{extension} of $A$.
The \textit{closure} $\overline{A}$ of $A$,
when it exists,
is the operator characterized by 
\begin{align*}
    \{(x, \overline{A} x)\ |\ x \in \mathrm{Dom}(\overline{A})\} = \overline{\{(x, A x) \ |\ x \in \mathrm{Dom}(A)\}}.
\end{align*}
Here,
the closure on the right-hand side is taken in $\mathcal{H} \times \mathcal{H}$.
Writing $A^*$ for the adjoint operator,
we call $A$ \textit{symmetric} if $A \subset A^*$,
and \textit{self-adjoint} if $A = A^*$.
A symmetric operator is called \textit{essentially self-adjoint} if $\overline{A}$ is self-adjoint.

For $S \subset \mathcal{H}$,
we write $S^{\perp} := \{y \in \mathcal{H} \ |\ \langle y,x \rangle_{\mathcal{H}} = 0 \text{ for all } x\in S\}$,
and set 
\begin{align*}
    \mathrm{ran}\,A := \{A x \ |\ x \in \mathrm{Dom}(A) \},\quad \mathrm{ker} \,A := \{x \in \mathrm{Dom}(A)\ |\ Ax = 0\}.
\end{align*}
If $A$ is self-adjoint,
then we have 
\begin{align}\label{eq:range}
    \overline{\mathrm{ran}\,A} = (\mathrm{ker}\,A)^{\perp},
\end{align}
where the left-hand side is the closure of $\mathrm{ran} \,A$ in $\mathcal{H}$.
\subsection{Functional analysis for manifolds with boundary}\label{subsec:functional-analysis-boundary}
We introduce the Neumann Laplacian through its energy form as in Bianchi--G\"{u}neysu--Setti \cite{bianchi2025neumann}.
The approximation results recalled below will be used to derive the lower bound for $\lambda_1^{\mathrm{Neu}}(M,g)$.
We denote the Riemannian volume measure of $g$ and $g_{\partial M}$ by $v_g$ and $v_{g,\partial M}$,
respectively.
The Laplacian is defined as follows:
\begin{align*}
    \Delta \varphi := \mathrm{div}(\nabla \varphi).
\end{align*}
For simplicity,
we write $L^p(M,v_g)$ as $L^p(M)$ for $p = 1,2$.
Also,
we write the inner product of $L^2(M)$ as $\langle \cdot, \cdot \rangle_{L^2}$,
and the norm of this inner product as $\|\cdot\|_2$.
Let $(N,g_N)$ be a Riemannian manifold without boundary.
For $p = 1,2$,
we denote the Sobolev space as follows (see Bianchi--G\"{u}neysu--Setti \cite{bianchi2025neumann}):
\begin{align*}
    W^{1,p}(N) := \{ \varphi \in L^p(N)\ |\ \nabla \varphi \in \Gamma_{L^p}(N,TN) \},
\end{align*}
where we refer to \cite{guidetti2017elliptic} for the definition of $\Gamma_{L^p}(N,TN)$.
For $\mathring{M} := M \backslash \partial M$,
we set $\mathrm{Dom} (\mathcal{E}_g) := W^{1,2}(\mathring{M})$,
and define 
\begin{align*}
    \mathcal{E}_g(\varphi,\psi) := \int_M g(\nabla \varphi, \nabla \psi)\,\d v_g
\end{align*}
for $\varphi, \psi \in \mathrm{Dom}(\mathcal{E}_g)$.
Let $L_g$ be the uniquely determined non-negative self-adjoint operator in $L^2(M)$ satisfying $\mathrm{Dom}(L_g) \subset W^{1,2}(\mathring{M})$ and
\begin{align*}
    \langle L_g \varphi, \psi\rangle_{L^2}  = \mathcal{E}_g(\varphi,\psi)
\end{align*}
for $\varphi \in \mathrm{Dom}(L_g)$ and $\psi \in W^{1,2}(\mathring{M})$.
\textcolor{black}{We call $L_g$ the \textit{Neumann Laplacian}.}
\begin{remark}\label{rem:kernel}
We first observe that 
\begin{align}\label{eq:lem:upper-5}
     \mathrm{ker}\, L_g = \left\{ \varphi \in L^2(M)\ |\ \varphi \text{ is constant }v_g\text{-almost everywhere} \right\}.
\end{align}
Indeed,
if $\varphi \in \mathrm{ker} \,L_g$,
then we have 
\begin{align*}
    \int_M |\nabla \varphi|^2\,\d v_g = \langle L_g \varphi ,\varphi \rangle_{L^2} = 0.
\end{align*}
Thus,
the weak gradient vanishes.
This implies that $\varphi$ is constant $v_g$-almost everywhere.
Conversely,
let $\varphi \in L^2(M)$ be constant $v_g$-almost everywhere.
We have $\varphi \in W^{1,2}(\mathring{M})$ and $\mathcal{E}_g(\varphi,\psi) = 0$ for every $\psi \in W^{1,2}(\mathring{M})$. 
It follows from the definition of $L_g$ that $\varphi \in \mathrm{Dom}(L_g)$ and $L_g \varphi = 0$.
This proves \eqref{eq:lem:upper-5}.
\end{remark}
Let $C_c^\infty(M)$ be the set of compactly supported functions in $C^\infty(M)$.
We set 
\begin{align*}
    \widehat{C}_c^\infty(M) := \{ \varphi \in C_c^\infty(M)\ |\ g(\nabla \varphi, \nu_{\partial M}) = 0 \},
\end{align*}
where $\nu_{\partial M}$ denotes the outward unit normal vector field on $\partial M$.
Regarding the density between functional spaces on manifolds with boundary,
we have the following assertion (see Bianchi--G\"{u}neysu--Setti \cite[Theorem 1.4]{bianchi2025neumann}):
\begin{theorem}[\cite{bianchi2025neumann}]\label{thm:dense}
Let $(M,g)$ be a complete Riemannian manifold with boundary.
Then $\widehat{C}_c^\infty(M)$ is dense in $W^{1,2}(\mathring{M})$.
\end{theorem}
\begin{remark}\label{rem:dense}
By integration by parts and Theorem \ref{thm:dense},
we have 
\begin{align*}
    \mathcal{E}_g (\varphi,\psi) = \langle -\Delta\varphi,\psi \rangle_{L^2}
\end{align*}
for $\varphi \in \widehat{C}_c^\infty(M)$ and $\psi \in W^{1,2}(\mathring{M})$.
Since $\Delta \varphi \in L^2(M)$,
the definition of $L_g$ yields 
\begin{align}\label{eq:rem:dense}
    \widehat{C}_c^\infty (M) \subset \mathrm{Dom}(L_g),\quad L_g \varphi = -\Delta \varphi.
\end{align}
\end{remark}
Regarding the approximation results,
we also have the following essential self-adjointness theorem (see Bianchi--G\"{u}neysu--Setti \cite[Theorem 2.1]{bianchi2025neumann}):
\begin{theorem}[\cite{bianchi2025neumann}]
Let $(M,g)$ be a complete Riemannian manifold with boundary.
Then $-\Delta|_{\widehat{C}_c^\infty(M)}$ is essentially self-adjoint.
In particular,
$L_g$ is the unique self-adjoint extension of this operator.
\end{theorem}
\begin{remark}\label{rem:graph}
This theorem implies 
\begin{align*}
    \overline{-\Delta|_{\widehat{C}_c^\infty(M)}} = L_g.
\end{align*}
Hence,
for every $\varphi \in \mathrm{Dom}(L_g)$,
there exists a sequence $(\varphi_j)$ in $\widehat{C}_c^\infty(M)$ such that 
\begin{align*}
    \varphi_j \longrightarrow \varphi,\quad L_g \varphi_j \longrightarrow L_g \varphi 
\end{align*}
in $L^2(M)$ (see also the proof of \cite[Theorem 2.1]{bianchi2025neumann} and \cite[Definition 1.5]{schmudgen2012unbounded}).
Together with \eqref{eq:range},
we see 
\begin{align}\label{eq:rem:graph}
    \overline{L_g \left( \widehat{C}_c^\infty (M) \right)} = \overline{\mathrm{ran}\,L_g} = (\mathrm{ker} \,L_g)^{\perp}.
\end{align}
\end{remark}
\subsection{Geometric inequalities}\label{subsec:geometric-inequality}
We next recall the geometric estimates used in the next section.
We define the second fundamental form and the mean curvature as follows:
\begin{align*}
    \mathrm{II}_{g,\partial M}(X,Y) := g(\nabla_X \nu_{\partial M}, Y),\quad H_{g,\partial M} := \mathrm{tr} \,\mathrm{II}_{g,\partial M}.
\end{align*}
Let $d_g$ denote the distance function induced by the Riemannian metric $g$ on $M$ (see e.g., Bianchi--G\"{u}neysu--Setti \cite[Section 1]{bianchi2025neumann} and Sakurai \cite[Section 2]{sakurai2017rigidity}).
For $t \geq 0$ and $p \in M$,
we set 
\begin{align*}
    B_t(p) := \{ x \in M \ | \ d_g(p,x) < t \}.
\end{align*}
We have the following Bishop--Gromov volume comparison (see Mascellani \cite[Proposition 3.4]{mascellani2017fourth}):
\begin{proposition}[\cite{mascellani2017fourth}]\label{prop:bishop-gromov}
Let $(M,g)$ be an $n$-dimensional complete Riemannian manifold with boundary.
We assume 
\begin{align*}
    \Ric_g \geq 0,\quad \mathrm{II}_{g,\partial M} \geq 0.
\end{align*}
Then for every $p \in M$,
the function
\begin{align*}
    t \longmapsto \frac{v_g(B_t(p))}{t^n}
\end{align*}
is non-increasing on $(0,\infty)$.
\end{proposition}
We set $\rho_{\partial M} (x) := d_g(x,\partial M)$.
We have the following inscribed radius estimate (see Kasue \cite[Theorem A]{kasue1983ricci} and Li \cite[Theorem 1.1]{li2014sharp}):
\begin{theorem}[\cite{kasue1983ricci,li2014sharp}]\label{thm:li}
Let $(M,g)$ be an $n$-dimensional complete Riemannian manifold with boundary.
For a constant $\kappa > 0$,
we assume 
\begin{align*}
    \Ric_g \geq 0,\quad H_{g,\partial M} \geq (n-1)\kappa.
\end{align*}
Then we have 
\begin{align*}
    \sup_{x\in M} \rho_{\partial M}(x) \leq \kappa^{-1}.
\end{align*}
\end{theorem}
Below,
we record the special case of Yan--Zhu \cite[Proposition A.3]{yan2026nonnegative},
which is sufficient for our purpose.
Indeed,
we have the following $L^1$-trace Poincar\'{e} inequality:
\begin{proposition}[\cite{yan2026nonnegative}]\label{prop:appendix}
Let $(M,g)$ be an $n$-dimensional complete Riemannian manifold with boundary.
We assume 
\begin{align*}
    \Ric_g \geq 0,\quad H_{g,\partial M} > 0.
\end{align*}
Let $\pi(x) \in \partial M$ be the nearest boundary point, 
defined for almost every $x\in M$.
For $f \in W^{1,1}(\mathring{M})$,
we have 
\begin{align*}
    \int_M |f|\,\d v_g \leq \frac{n-1}{n} \int_{\partial M} \frac{|f|}{H_{g,\partial M}} \,\d v_{g,\partial M} + \frac{1}{n} \int_M \left( \frac{n-1}{H_{g,\partial M}\circ \pi} - \rho_{\partial M}\right) |\nabla f|\,\d v_g,
\end{align*}
where the boundary values of $f$ are understood as traces.
\end{proposition}
\section{Proof of the main theorem}\label{sec:proof}
In this section,
we prove Theorem \ref{thm:yan-zhu} by comparing lower and upper bounds for the infimum of the Neumann Rayleigh quotient.
\subsection{Lower bound}\label{subsec:lower-bound}
We define 
\begin{align}\label{eq:def-lambda}
    \lambda_1^{\mathrm{Neu}}(M,g) := \inf_{\varphi} \frac{\mathcal{E}_g(\varphi, \varphi)}{\|\varphi\|_2^2},
\end{align}
where the infimum is taken over all $0 \neq \varphi \in W^{1,2}(\mathring{M}) \cap (\mathrm{ker}\,L_g)^{\perp}$.
\begin{remark}\label{rem:lambda}
We consider the case of finite volume $v_g(M) < \infty$.
We set $\mu := v_g(M)^{-1} v_g$.
By Remark \ref{rem:kernel},
the condition $\varphi \in (\mathrm{ker}\,L_g)^{\perp}$ is equivalent to $\int_M \varphi \,\d \mu = 0$.
Subtracting the mean of $\varphi$ gives 
\begin{align*}
    \lambda_1^{\mathrm{Neu}}(M,g) = \inf_{\substack{\varphi \in W^{1,2}(\mathring{M})\\
    \varphi \text{ is not constant a.e.}}} \frac{\int_M |\nabla \varphi|^2 \,\d \mu}{\int_M \varphi^2 \,\d \mu - \left( \int_M \varphi \,\d \mu\right)^2}.
\end{align*}
This coincides with the unweighted case of the variational quantities in Wang \cite{wang1999existence}.
Indeed,
although Wang \cite{wang1999existence} takes the infimum over $C^1(M) \cap L^2(M,\mu)$,
replacing it with $W^{1,2}(\mathring{M})$ does not change the infimum since we have Theorem \ref{thm:dense}.
\end{remark}
We first derive an estimate for functions in $\widehat{C}_c^\infty(M)$. 
We then use the approximation property from Subsection \ref{subsec:functional-analysis-boundary} to obtain a lower bound of $\lambda_1^{\mathrm{Neu}}(M,g)$.
We first show the following inequality:
\begin{lemma}\label{lem:eigenvalue-estimate}
Let $(M,g)$ be an $n$-dimensional complete Riemannian manifold with boundary.
For $\kappa > 0$,
we assume
\begin{align}\label{eq:lem:curvature-assumptions}
    \Ric_g \geq 0,\quad \mathrm{II}_{g,\partial M} \geq \kappa g_{\partial M}.
\end{align}
Then for $\varphi \in \widehat{C}_c^\infty(M)$,
we have 
\begin{align*}
    \int_M (\Delta \varphi)^2 \,\d v_g \geq  (n-1)\kappa^2 \int_M |\nabla \varphi|^2 \,\d v_g.
\end{align*}
\end{lemma}
\begin{proof}
We have the following Reilly formula (see Reilly \cite{reilly1977applications}):
\begin{align*}
    \int_M (\Delta \varphi)^2 \,\d v_g &= \int_M \left( |\Hess \varphi|^2 + \Ric_g(\nabla \varphi, \nabla \varphi) \right)\,\d v_g + \int_{\partial M} H_{g,\partial M} g(\nabla \varphi, \nu_{\partial M})^2 \,\d v_{g,\partial M}\\
    &\,\,\quad + 2 \int_{\partial M} g(\nabla \varphi, \nu_{\partial M}) \Delta_{\partial M}\,\psi\,\d v_{g,\partial M} + \int_{\partial M} \mathrm{II}_{g,\partial M} (\nabla_{\partial M}\,\psi, \nabla_{\partial M}\,\psi)\,\d v_{g,\partial M},
\end{align*}
where $\Delta_{\partial M}$ and $\nabla_{\partial M}$ are Laplacian and gradient with respect to $g_{\partial M}$,
and $\psi := \varphi|_{\partial M}$.
Together with $g(\nabla \varphi, \nu_{\partial M}) = 0$,
$\nabla \varphi = \nabla_{\partial M} \,\psi$ and the assumptions on the curvatures \eqref{eq:lem:curvature-assumptions},
we have 
\begin{align}\label{eq:lem:eigenvalue-estimate-1}
    \int_M (\Delta \varphi)^2 \,\d v_g \geq \int_M |\Hess \varphi|^2 \,\d v_g + \kappa \int_{\partial M} |\nabla \varphi|^2 \,\d v_{g,\partial M}.
\end{align}
Also,
we note that Theorem \ref{thm:li} gives 
\begin{align}\label{eq:lem:eigenvalue-estimate-1-5}
    1 - \kappa \rho_{\partial M} \geq 0.
\end{align}
Applying Proposition \ref{prop:appendix} with $f = |\nabla \varphi|^2$,
and using \eqref{eq:lem:eigenvalue-estimate-1-5} and
\begin{align*}
    |\nabla f| = 2 | \Hess \varphi (\nabla \varphi, \cdot )| \leq 2 |\Hess \varphi| |\nabla \varphi|,\quad H_{g,\partial M} \geq (n-1)\kappa,
\end{align*}
we obtain 
\begin{align*}
    n\kappa^2 \int_M |\nabla \varphi|^2\,\d v_g \leq \kappa \int_{\partial M} |\nabla \varphi|^2 \,\d v_{g,\partial M} + 2 \kappa \int_M (1-\kappa \rho_{\partial M}) |\Hess \varphi| |\nabla \varphi|\,\d v_g.
\end{align*}
Together with Young's inequality and \eqref{eq:lem:eigenvalue-estimate-1},
we see 
\begin{align*}
    n\kappa^2 \int_M |\nabla \varphi|^2\,\d v_g &\leq \kappa \int_{\partial M} |\nabla \varphi|^2 \,\d v_{g,\partial M}  + \int_M |\Hess \varphi|^2 \,\d v_g + \kappa^2 \int_M (1-\kappa \rho_{\partial M})^2 |\nabla \varphi|^2\,\d v_g\\
    &\leq \int_M (\Delta\varphi)^2 \,\d v_g + \kappa^2 \int_M (1-\kappa \rho_{\partial M})^2 |\nabla \varphi|^2 \,\d v_g.
\end{align*}
Combining this with 
\begin{align*}
    \kappa^2 \left\{ n - (1-\kappa \rho_{\partial M})^2 \right\} = (n-1) \kappa^2 + \kappa^3 \rho_{\partial M} (2 - \kappa \rho_{\partial M}),
\end{align*}
we have 
\begin{align*}
    \int_M (\Delta \varphi)^2 \,\d v_g - (n-1)\kappa^2 \int_M |\nabla \varphi|^2\,\d v_g \geq \kappa^3 \int_M \rho_{\partial M} (2 - \kappa \rho_{\partial M}) |\nabla \varphi|^2 \,\d v_g.
\end{align*}
The right-hand side is non-negative by \eqref{eq:lem:eigenvalue-estimate-1-5}.
We conclude the proof.
\end{proof}
Combining this with the approximation argument in Subsection \ref{subsec:functional-analysis-boundary},
we obtain the lower bound for $\lambda_1^{\mathrm{Neu}}(M,g)$ as follows:
\begin{lemma}\label{lem:poincare}
Let $(M,g)$ be an $n$-dimensional complete Riemannian manifold with boundary.
For $\kappa > 0$,
we assume
\begin{align*}
    \Ric_g \geq 0,\quad \mathrm{II}_{g,\partial M} \geq \kappa g_{\partial M}.
\end{align*}
Then we have
\begin{align}\label{eq:lem:poincare-1}
    \lambda_1^{\mathrm{Neu}}(M,g) \geq (n-1) \kappa^2.
\end{align}
\end{lemma}
\begin{proof}
We consider nonzero $\varphi \in W^{1,2}(\mathring{M}) \cap (\mathrm{ker}\,L_g)^{\perp}$.
It follows from \eqref{eq:rem:graph} that
\begin{align*}
    \varphi \in \overline{L_g \left( \widehat{C}_c^\infty (M) \right)}
\end{align*}
Hence,
there exists a sequence $(f_j)$ in $\widehat{C}_c^\infty(M)$ such that 
\begin{align*}
    L_g f_j \longrightarrow \varphi
\end{align*}
in $L^2(M)$ when $j \rightarrow \infty$.
It follows from the definition of $L_g$,
the Cauchy--Schwarz inequality and Lemma \ref{lem:eigenvalue-estimate} that 
\begin{align*}
    \left| \langle L_g f_j, \varphi \rangle_{L^2} \right|^2 = \left| \mathcal{E}_g(f_j,\varphi) \right|^2 \leq \mathcal{E}_g (f_j, f_j)\, \mathcal{E}_g(\varphi, \varphi) \leq \frac{\|L_g f_j\|^2_2}{(n-1)\kappa^2} \mathcal{E}_g(\varphi, \varphi).
\end{align*}
Passing to the limit as $j \rightarrow \infty$,
we have 
\begin{align*}
    \|\varphi\|_2^4 \leq \frac{\|\varphi\|_2^2}{(n-1)\kappa^2} \mathcal{E}_g(\varphi, \varphi).
\end{align*}
This gives the desired inequality.
%
%
\end{proof}
\subsection{Upper bound}\label{subsec:upper-bound}
We have the following upper bound for the Neumann Rayleigh quotient:
\begin{lemma}\label{lem:upper}
Let $(M,g)$ be an $n$-dimensional complete Riemannian manifold with boundary.
We assume
\begin{align*}
    \Ric_g \geq 0,\quad \mathrm{II}_{g,\partial M} \geq 0.
\end{align*}
For $p,q \in \mathring{M}$ with $p \neq q$,
there exists a nonzero function $\varphi \in W^{1,2}(\mathring{M}) \cap (\mathrm{ker}\,L_g)^{\perp}$ such that 
\begin{align*}
     \frac{\mathcal{E}_g(\varphi,\varphi)}{\|\varphi\|_2^2} \leq \frac{2(n + 1)(n + 2)}{d_g(p,q)^2}.
\end{align*}
\end{lemma}
\begin{proof}
For brevity,
we set $R := d_g(p,q)/2$.
For $y \in \{p,q\}$ and $x \in M$,
we set 
\begin{align*}
    r_y (x) := d_g(y,x),\quad
    \varphi_y(x) := \left( R - r_y(x)\right)_{+},
\end{align*}
where $(s)_+ :=\max\{s,0\}$ for $s \in \mathbb{R}$.
Let $\mathbf{1}_S$ denote the characteristic function of a set $S$. 
Since $\nabla \varphi_y(x) = - \mathbf{1}_{B_R(y)}(x)\,\nabla r_y (x)$ for $v_g$-a.e. $x \in M$,
we have 
\begin{align}\label{eq:lem:upper-0}
    |\nabla \varphi_y(x)| \leq \mathbf{1}_{B_R(y)}(x)
\end{align}
for $v_g$-a.e. $x \in M$.
This implies 
\begin{align}\label{eq:lem:upper-1}
    \mathcal{E}_g(\varphi_y, \varphi_y) \leq v_g(B_R(y)).
\end{align}
Also,
by Tonelli's theorem and Proposition \ref{prop:bishop-gromov}, 
we have 
\begin{align*}
    \|\varphi_y\|_2^2 &= \int_M \int_0^R 2(R -t) \mathbf{1}_{B_{t}(y)} \,\d t \d v_g(x)\nonumber\\
    &= 2\int_0^R (R - t) v_g(B_t(y))\,\d t\nonumber\\
    &\geq \frac{2 v_g(B_R(y))}{R^n} \int_0^R (R - t) t^n \,\d t\nonumber\\
    &= \frac{2R^2 v_g(B_R(y))}{(n + 1)(n + 2)}.\nonumber
\end{align*}
Combining this with \eqref{eq:lem:upper-1},
we have 
\begin{align}\label{eq:lem:upper-3}
    \mathcal{E}_g (\varphi_y, \varphi_y) \leq \frac{2(n + 1)(n + 2)}{d_g(p,q)^2} \|\varphi_y\|_2^2.
\end{align}
Since $2 R = d_g(p,q)$,
we see 
\begin{align*}
    \{x \in M\ |\ \varphi_p(x) > 0\} \cap \{ x\in M\ |\ \varphi_q(x) > 0 \} = \emptyset.
\end{align*}
This implies 
\begin{align}\label{eq:lem:upper-4}
    \langle \varphi_p,\varphi_q \rangle_{L^2} = 0,\quad \mathcal{E}_g(\varphi_p,\varphi_q) = 0.
\end{align}
We set 
\begin{align*}
    m_y := \int_M \varphi_y \,\d v_g,\quad \varphi := \frac{\varphi_p}{m_p} - \frac{\varphi_q}{m_q}.
\end{align*}
We note that $\varphi \neq 0$.
By \eqref{eq:lem:upper-0},
we have $\varphi \in W^{1,2}(\mathring{M})$.
Next,
we show $\varphi \in (\mathrm{ker}\,L_g)^{\perp}$.
Let $\psi \in \mathrm{ker}\,L_g$.
By Remark \ref{rem:kernel},
there exists $c \in \mathbb{R}$ such that $\psi(x) = c$ for $v_g$-a.e. $x \in M$.
Together with the definition of $\varphi$,
we have 
\begin{align*}
    \langle \psi, \varphi \rangle_{L^2} = c \int_M \varphi\,\d v_g = 0.
\end{align*}
Therefore,
we have $\varphi \in (\mathrm{ker}\,L_g)^{\perp}$.
Lastly,
it follows from \eqref{eq:lem:upper-3} and \eqref{eq:lem:upper-4} that we have 
\begin{align*}
    \mathcal{E}_g(\varphi, \varphi) \leq \frac{2(n + 1)(n + 2)}{d_g(p,q)^2}\left( \frac{\|\varphi_p\|_2^2}{m_p^2} + \frac{\|\varphi_q\|_2^2}{m_q^2} \right) = \frac{2(n + 1)(n + 2)}{d_g(p,q)^2} \| \varphi\|_2^2.
\end{align*}
We conclude the proof.
\end{proof}
We are now in a position to give a proof of the main theorem:
\begin{proof}[Proof of Theorem \ref{thm:yan-zhu}]
We fix distinct $p,q \in \mathring{M}$.
Lemmas \ref{lem:poincare} and \ref{lem:upper} give
\begin{align*}
    (n-1)\kappa^2 
    \leq \lambda_1^{\mathrm{Neu}}(M,g)
    \leq \frac{2(n + 1)(n + 2)}{d_g(p,q)^2}.
\end{align*}
This yields
\begin{align}\label{eq:thm:main-1}
    d_g(p,q) \leq \frac{1}{\kappa} \sqrt{\frac{2(n + 1)(n + 2)}{n-1}}.
\end{align}
This implies that the diameter of $M$ is bounded above by the right-hand side of \eqref{eq:thm:main-1}.
Since $(M,d_g)$ is a complete and locally compact length space,
the Hopf--Rinow theorem implies that every closed bounded subset of $M$ is compact.
Therefore,
the diameter bound yields that $M$ is compact.
\end{proof}
\begin{remark}
We note that the finiteness of the fundamental group $\pi_1(M)$ also follows by the same argument as in Yan--Zhu \cite{yan2026nonnegative}.
\end{remark}
\subsection*{Acknowledgements}
The author expresses sincere gratitude to Professor Chao Xia for valuable comments and encouragement.
The author also thanks Professor Xingyu Zhu for drawing the author's attention to Ref.~\cite{yan2026nonnegative}.
Part of this work was carried out during the author's stay at the Tianyuan Mathematical Center in Southeast China at Xiamen University.
The author is grateful to the center for its hospitality,
and acknowledges financial support for this stay from the UTokyo Global Activity Support Program for Young Researchers.
This work was supported by JSPS KAKENHI Grant Number JP25KJ0271.
\subsection*{Use of AI}
The author prepared this manuscript with the aid of ChatGPT.
The author takes full responsibility for the content of this paper.

\bibliographystyle{amsplain}
\bibliography{ref}
\end{document}